\documentstyle[11pt,psfig]{article} 
\input{psfig} 
 \newtheorem{theorem}{Theorem}[section] 
 \newtheorem{lemma}[theorem]{Lemma} 
 \newtheorem{corollary}[theorem]{Corollary}

\newtheorem{definition}[theorem]{Definition}
\newtheorem{proposition}[theorem]{Proposition}

\newenvironment{proof}{{\it Proof:\/}}{$\Box$\vskip 0.08in}

\begin{document} 
\thispagestyle{empty} 
\ 
\vspace{0.5in} 
 \begin{center} 
 {\LARGE\bf Kauffman bracket skein module of a connected sum of 3-manifolds}
\end{center} 
\vspace*{0.5in} 
 \begin{center} 
                      J\'ozef H.~Przytycki 
\end{center} 
\ \\
\vspace*{0.5in} 
{\bf Abstract.}\\
We show that for the Kauffman bracket skein module over the field of rational
functions in variable $A$, the module of a connected sum of 3-manifolds is
the tensor product of modules of the individual manifolds.
\section{The main theorem.}\label{1}
We recall in this section the definition of the Kauffman bracket skein module
(KBSM) and formulate the main result of the paper.
\begin{definition}[\cite{P-1,H-P-1}]\label{1.1}\ \\ 
Let $M$ be an oriented 3-manifold, $R$ a commutative ring with identity and $A$ 
its invertible element. Let ${\cal L}_{fr}$ be the set of unoriented 
framed links in $M$, including the empty link ${\emptyset}$. 
Let $S$ be 
the submodule of $R{\cal L}_{fr}$ generated by skein expressions 
$L_+- AL_0 - A^{-1}L_{\infty}$, where the triple 
$L_+, L_0, L_{\infty}$ is shown in Fig.1.1, and $L \sqcup T_1 +(A^2 + 
A^{-2})L$, where $T_1$ denotes the trivial framed knot. 
We define the Kauffman bracket skein module (KBSM), 
denoted by\footnote{
The standard notation for the KBSM is ${\cal S}_{2,\infty}(M;R,A)$,
\cite{P-1,H-P-1}, but in this paper we do not discuss 
skein modules other than KBSM so for simplicity we 
drop $(2,\infty)$ from the notation.}
 ${\cal S}(M;R,A)$, 
as the quotient ${\cal S}(M;R,A)= R{\cal L}_{fr}/S$. 
\end{definition} 
\ \\ 
\centerline{\psfig{figure=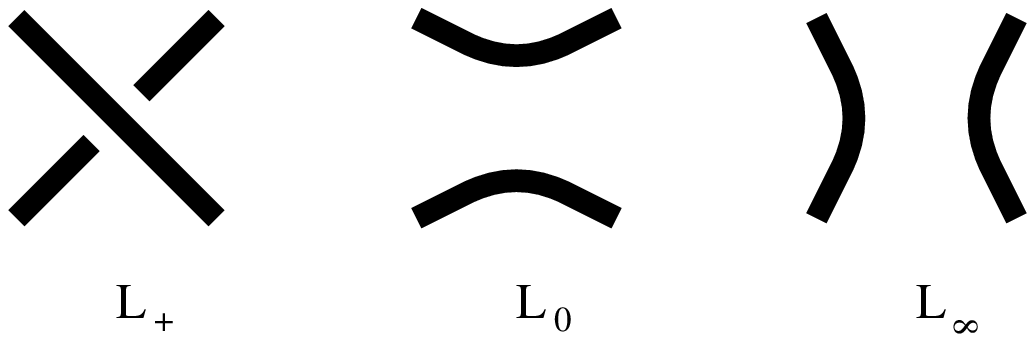,height=1.6cm}} 
\begin{center} 
Fig. 1.1. 
\end{center} 
Notice that $L^{(1)}=-A^3L$ in ${\cal S}(M;R,A)$, where 
$L^{(1)}$ denotes a link obtained from  $L$ by one positive twist of 
the framing of $L$.  We call this the framing relation. 
For the sake of shortness of notation, we will often drop $(R,A)$ from
${\cal S}(M;R,A)$ and write simply ${\cal S}(M)$, as long as it 
unambiguous. 

\begin{theorem}\label{1.2} 
Assume that $(A^k-1)$ is invertible in $R$ for any $k>0$. Then
$${\cal S}(M_1\#M_2)={\cal S}(M_1)\otimes {\cal S}(M_2),$$
where $M_1\#M_2$ denotes the connected sum of compact 
3-manifolds $M_1$ and $M_2$.
\end{theorem}

In particular we have:
\begin{corollary}
If $R$ is a field of rational functions in variable $A$, ${\cal F}(A)$,
or $R$ is a field of complex numbers, $C$, and $A$ is not a root
of unity, then ${\cal S}(M_1\#M_2)={\cal S}(M_1)\otimes {\cal S}(M_2)$. 
\end{corollary}

\section{Basic properties of skein modules}\label{2}
Below we list some elementary properties of KBSM (which also hold for
 other skein modules), \cite{P-1} (compare also \cite{H-P-1,P-3}).
\begin{proposition}\label{2.1}
\begin{enumerate}
\item [(1)] An orientation preserving embedding of 3-manifolds 
$i: M \to N$ yields the homomorphism of skein modules 
$i_*: {\cal S}(M) \to {\cal S}(N)$. 
The above correspondence leads to a functor from the category of 
3-manifolds and orientation preserving embeddings (up to ambient isotopy) 
to the category of $R$-modules (with a specified invertible element $A \in R$).
\item [(2)] 
\begin{enumerate}
\item [(i)] 
If $N$ is obtained from $M$ by adding a 3-handle to it (i.e. capping off
a hole so that $M=N\# D^3$), and $i: M \to N$ is the associated embedding, 
then $i_*:{\cal S}(M) \to {\cal S}(N)$
is an isomorphism.
\item [(ii)] If $N$ is obtained from $M$ by adding a 2-handle to it,
and $i: M \to N$ is the associated embedding, 
then $i_*:{\cal S}(M) \to {\cal S}(N)$    
is an epimorphism.
\end{enumerate}
\item [(3)] If $M_1 \sqcup M_2$ is the disjoint sum of 3-manifolds $M_1$
and $M_2$ then 
$${\cal S}(M_1\sqcup M_2)= {\cal S}(M_1) \otimes {\cal S}(M_2).$$  
\item [(4)]  
(Universal Coefficient Property)\\ 
Let $r: R \to R'$ be a homomorphism of rings (commutative with 1). 
We can think of $R'$ as an $R$ module. Then the identity map on 
${\cal L}_{fr}$ induces the isomorphism of $R'$ (and $R$) modules: 
$$ \bar r: {\cal S}(M;R,A)\otimes_R R' \to
{\cal S}(M;R',r(A))  .$$
\item [(5)] If $F$ is a surface, then the KBSM 
${\cal S}(F\times I)$ is a free $R$-module with  
basis $B(F)$ consisting of links on $F$, up to ambient isotopy of $F$,
  without contractible 
components (but including the empty link).     
\end{enumerate}
\end{proposition}
Results in Proposition 2.1 are well known; compare \cite{P-1,H-P-1,P-S,P-3}. 
We clarify some points of them below:
\begin{enumerate}
\item [(1)] If $i:M \to N$ is an orientation reversing
embedding then $i_*$ is a $Z$-homomorphism and $i_*(Aw)=A^{-1}i_*(w)$.
\item [(2)]
\begin{enumerate}
\item [(i)] It holds because the co-core of a 3-handle is $0$-dimensional.
\item [(ii)] It holds because the co-core of a 2-handle is $1$-dimensional.
\end{enumerate}
\item [(3)] This is a consequence of the well known property of short exact 
sequences, \cite{Bl}:\\
If $ 0 \to A' \to A \to A'' \to 0$ and $ 0 \to B' \to B \to B'' \to 0$ are
short exact sequences of $R$-modules then
$ 0 \to A'\otimes B + A\otimes B' \to A\otimes B \to A''\otimes B'' \to 0$ 
is a short exact sequence. 
\item [(4)] This important fact follows easily from right exactness of
the tensor functor (applied to a short exact sequence) and from 
the ``five lemma" (see for example \cite{C-E}).
\item [(5)] This applies in particular to a handlebody, because
$H_n = P_n \times I$, where $H_n$ is a handlebody of genus $n$ and $P_n$ is
a disc with $n$ holes.
\end{enumerate}

\section{Outline of the proof of the main theorem.}\label{3}
\begin{enumerate}
\item[1.]
Any compact 3-dimensional manifold can be obtained from a handlebody
by adding 2 and 3-handles to it. The KBSM of a handlebody is a well
understood free module (Prop. 2.1(5)), adding a 3-handle  does not change
the module (Prop. 2.1(2)(i)) and adding a 2-handle gives new relations 
to the skein module, but no new generators (Prop. 2.1(2)(ii)).

\item[2.] 
If a 2-handle is added, all new relations are obtained by
sliding links along the 2-handle (Lemma 4.1). 
There is usually an infinite collection of
relations, but in the case of a 2-handle added along a meridian 
curve (creating $S^2$), we prove that the relations form a ``controllable" 
sequence, which, over the field ${\cal F}(A)$,
allows us to reduce all curves cutting the sphere, but not more. 

\item[3.] The embedding $j: M_1\# D^3 \sqcup M_2\# D^3  \to M_1\#M_2$ 
yields an epimorphism of the KBSM. To see that every link in 
${\cal S}(M_1\#M_2)$
is in the image, it suffices to consider relations given by very simple 
slidings (Fig. 6.1), using the ``second side" of $S^2$. 
We use also the fact that $A^k-1$ is invertible in $R$.  
$\sqcup$ denotes the disjoint sum. The connected sum $M\# D^3$ is a
manifold obtained from $M$ by cutting off a hole in $M$. In particular,
by Proposition 2.1(2)(i), the skein modules of $M$ and $M\# D^3$
coincide.

\item[4.]
We will start the proof of Theorem 1.2, by considering the case of $M_1$ and
$M_2$ being handlebodies; $M_1=H_n, M_2=H_m$. 
$H_n\# H_m$ is equal to $H_{n+m}$ with a 2-handle added along
the boundary of the meridian disc separating $H_n$ from $H_m$. We show  
that the embedding $i: H_n \sqcup H_m \to H_n\# H_m$ yields an isomorphism 
of skein modules assuming that $A^k-1$ is invertible in $R$.
For this we show that all sliding relations are generated by slidings of
Fig. 5.1.
\item[5.] 
We generalize 4. by considering any $M_1\#M_2$ and observing that 
$M_1\#M_2$ is obtained from $H_n\# H_m$ by adding 2-handles to $H_n$ or
$H_m$.
\item[6.]
In steps 4 and 5 we have to show that even a very complicated sliding
(say, a curve in $H_n$ first being pushed to $H_m$ and back several times
and only then slid), can be reduced to slidings of Fig. 5.1 or
slidings taking place totally in $H_n$ or $H_m$ (compare Lemma 6.1).
\end{enumerate}

\section{Handle sliding lemma.}\label{4}
\begin{lemma}\ \\
Let $(M,\partial M)$ be a 3-manifold with the boundary $\partial M$,
and let $\gamma $ be a simple closed curve on the boundary. Let $N=M_{\gamma}$
be the 3-manifold obtained from $M$ by adding
a 2-handle along $\gamma $. Furthermore let ${\cal L}_{fr}^{gen}$ be a set
of framed links in $M$ generating ${\cal S}(M)$.\\
Then ${\cal S}(N) = {\cal S}(M)/J$, where
$J$ is the submodule of ${\cal S}(M)$ generated by 
expressions $L- sl_{\gamma}(L)$, where $L\in {\cal L}_{fr}^{gen}$ and
$sl_{\gamma}(L)$ is obtained from $L$ by sliding it along ${\gamma}$
(i.e. handle sliding).
\end{lemma}
\begin{proof}
Let $S^1\times [-1,1]$ be a tubular neighborhood of $\gamma$ in $\partial M$
($\gamma = S^1\times \{0\}$), and consider a 2-handle added along 
$\gamma$, that is $D^2\times D^1$ and a homeomorphism $\phi:
\partial D^2\times D^1 \to S^1\times [-1,1]$. Then $N= M{\bigcup}_{\phi}
D^2\times D^1$. Let $f: M\to N$ be a natural embedding, then $f_*:
{\cal S}(M) \to {\cal S}(N)$ is an 
epimorphism, Prop. 2.1(2)(ii) (any link in $N$ can be pushed (ambient isotoped) 
to $M$). Furthermore any skein relation can be performed in $M$. 
The only difference between KBSM of $M$ and $N$ lies in the fact that some 
nonequivalent links in $M$ can be equivalent in $N$; 
the difference lies exactly in the possibility
of sliding a link in $M$ along the added 2-handle (that is $L$ is moving
from one side of the co-core of the 2-handle to another); compare Fig.4.1. 
The proof of Lemma 4.1 is completed.

\end{proof}
\ \\ 
\centerline{\psfig{figure=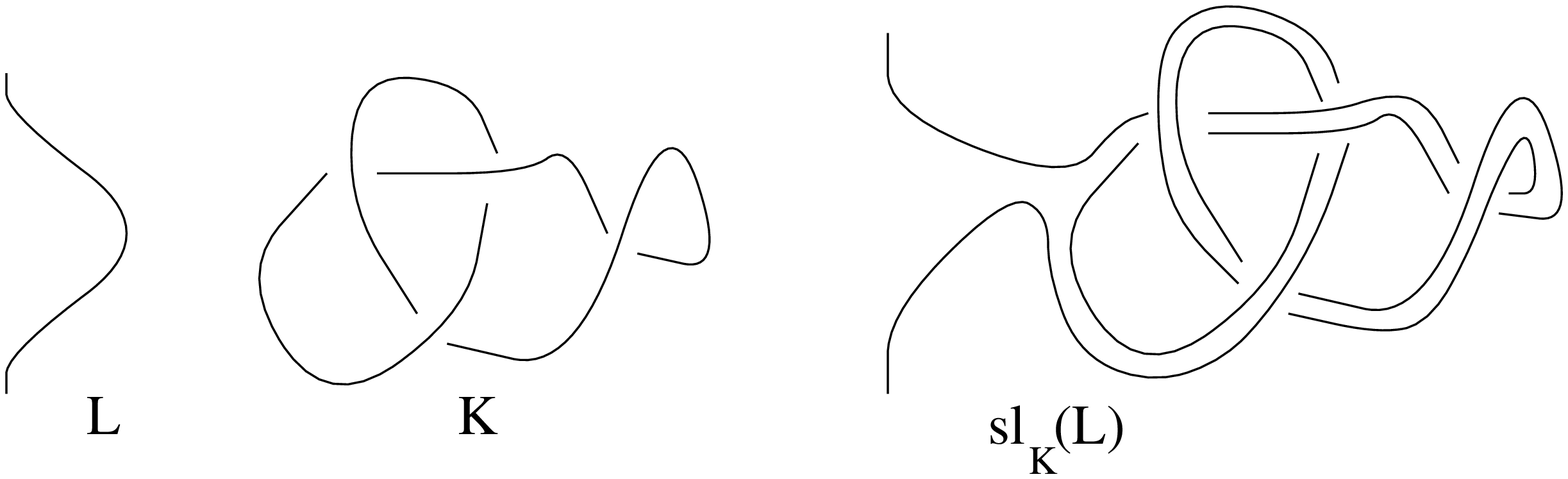,height=2.2cm}} 
\begin{center} 
Fig. 4.1. 
\end{center} 

Lemma 4.1 allows us to write an (infinite) presentation of the Kauffman bracket
skein module of any compact 3-manifold, using Heegaard decomposition and
knowledge of the module for handlebodies; Proposition 2.1(5).
This general presentation is not satisfactory and in some cases we can write
a simpler presentation.

\section{Epimorphism}\label{5}
We show in this section that for a connected sum $M_1\# M_2$, and $(A^k-1)$
invertible in $R$, ${\cal S}(M_1\# M_2)$ is generated by links,
components of which are in $M_1$ or $M_2$.\\

The above fact follows from the slightly more general proposition which
we prove below.
\begin{proposition}
Let $M$ be an oriented 3-manifold, $D$ a meridian disk in $M$, that is 
a properly embedded $2$-disk in $M$, and $\gamma = \partial D$.
If $R$ is a ring with $(A^k-1)$ invertible for any $k$ then the embedding
$j: (M- D) \to M_{\gamma}$, where $M_{\gamma}$ is obtained from
$M$ by adding a 2-handle along $\gamma$, induces an epimorphism
of the KBSM, $j_*: {\cal S}(M - D) \to {\cal S}(M_{\gamma})$.
\end{proposition}
\begin{proof}
The regular neighborhood, $V_D=[-1,1]\times D$, of $D$ in $M$ can be 
projected into 2-disk $D_p=[-1,1]\times [0,1]$ (then $V_D= D_p \times[0,1]$), 
and  we use $D_p$ to present link diagrams, compare Fig.5.1.
In ${\cal S}(M_{\gamma})$ one has sliding relations 
described in Fig. 5.1 (with blackboard framing).
These relations can be written as $p(z_k)= (-A^2-A^{-2})z_k$, 
where $z_k$ is a link in $M$ in general position with $D$ 
and cutting it $k$ times;
Fig.5.1.  After simplifying the formula, using the Kauffman bracket skein
relations, one gets: $p(z_k)= (-A^{2k+2}-A^{-2k-2})z_k +
\Sigma_{i=0}^{k-2}w_i(A)z_i$ and finally
$$(A^{2k+4}-1)(A^{2k}-1)z_k = A^{-2k-2}\Sigma_{i=0}^{k-2}w_i(A)z_i$$
that is $(A^{2k+4}-1)(A^{2k}-1)z_k$
is a linear combination of links with a smaller than $k$ intersection number
with the 2-sphere $D_{\gamma}$. For $(A^{2k+4}-1)(A^{2k}-1)$
invertible in $R$, one can eliminate $z_k$ from the set of generators.
Using induction, one can eliminate all elements of
${\cal S}(M_{\gamma})$ which cut
the 2-sphere $D_{\gamma}$ non-trivially. Thus $j_*$ is an epimorphism.
\end{proof}
\ \\
\centerline{\psfig{figure=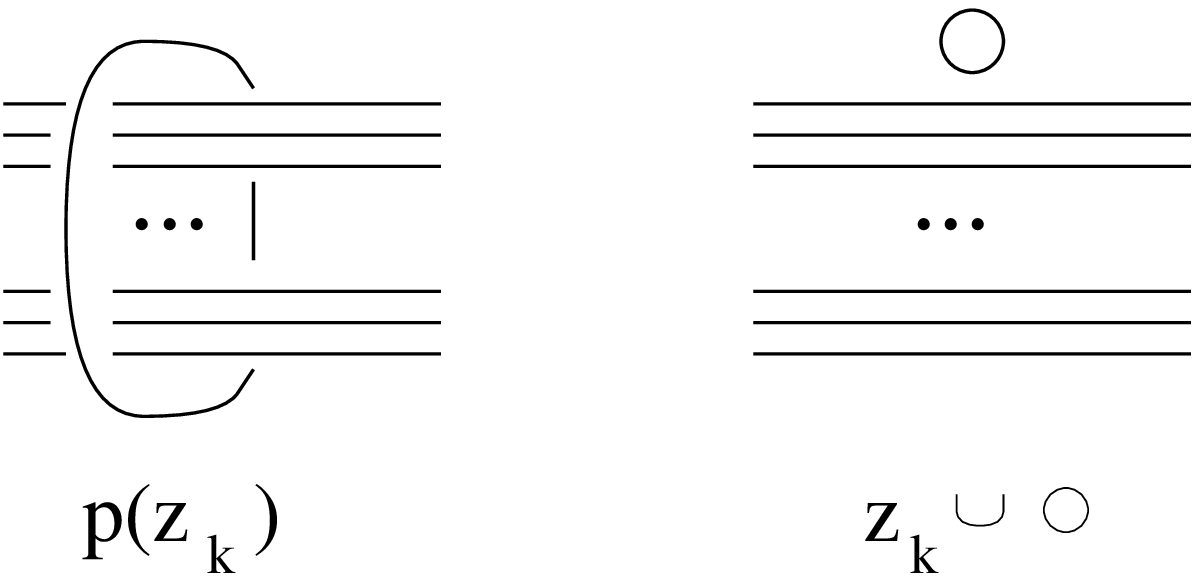,height=2.8cm}}
\begin{center}
Fig. 5.1.
\end{center}

\begin{corollary}\label{5.2}
If $R$ is a ring with $(A^k-1)$ invertible for any $k>0$
then the embedding $j: M_1\# D^3 \sqcup  M_2\# D^3 \to M_1\#M_2$ 
induces an epimorphism of the KBSM: 
$j_*: {\cal S}(M_1\# D^3 \sqcup  M_2\# D^3) \to {\cal S}(M_1\#M_2)$.
\end{corollary}

\section{Proof of the main theorem}\label{6}
We start the proof by showing that handle sliding described in Fig. 5.1
generates all handle sliding relations as long as $(A^k-1)$ is 
invertible for any $k>0$. This allows us to prove the main theorem for 
handlebodies as we know the basis of the KBSM in this case so we are
able to choose only those handle slidings which reduce those links from
the basis which cut $S^2$ (from the connected sum) non-trivially. Finally
we prove the main theorem for any compact 3-manifold using the fact that
such a manifold is obtained from a handlebody by adding 2 and 3-handles.

We say that a handle sliding $sl_{\gamma}$ of a link $L$ in $M$ along 
${\gamma}\in \partial M$ has support in a submanifold $V$ of $M$ if 
${\gamma}\in V$ and $L$ and $sl_{\gamma}(L)$ are identical outside $V$
\footnote{Here we consider the concrete realization of links not
ambient isotopy class. To omit confusion, we often will write
$D_L$ for a representative of the ambient isotopy class of a link $L$.}.
\begin{lemma}\label{6.1}
Let $D$ be a meridian disk in $M$, ${\gamma}= \partial D$, and
let $V_D = [-1,1] \times D$ be a regular neighborhood of $D$ in $M$.
\begin{enumerate}
\item[(a)]
If a link $L$ is disjoint from the disks $\{-1,1\} \times D$
and the sliding $sl_{\gamma}$ has support in $(-1,1)\times D$,
then the relation $L= sl_{\gamma}(L)$ holds in ${\cal S}(M)$.
\item[(b)]
Let $M_0$ be a component of $M$ which contains $D$ and assume that
$L \cap M_0$ is a trivial knot, then the sliding relation
$L = sl_{\gamma}(L)$ holds in ${\cal S}(M)$ for any sliding of $L$.
\item[(c)]
Let $L^{gen}_{fr}$ be a set of links generating ${\cal S}(M)$ and for
each $L$ in $L^{gen}_{fr}$ choose a representative $D_L$ embedded
in $M$ in such a way that $D_L$ cuts the meridian disk $D$ in the minimal
number of points in its ambient isotopy class. Let $I_0$ be a submodule of
${\cal S}(M)$ generated by sliding relations of Fig. 5.1, for $z_k = D_L$,
$L \in L^{gen}_{fr}$. In other words $I_0$ is generated 
by expressions $p(D_L) + (A^2 + A^{-2})D_L$.
Then for any representative, ${\tilde D}_L$, of a link $L$ in $M$
a sliding of Fig. 5.1 preserves the element of ${\cal S}(M)/I_0$ 
(i.e.  $(p(L) + (A^2 + A^{-2})L) \in I_0$).
\item[(d)]
If $R$ is a ring with $(A^k-1)$ invertible for any $k$  and $I_0$ is defined
as in (c) then any sliding relation holds in ${\cal S}(M)/I_0$.
\end{enumerate}
\end{lemma}
\begin{proof}
\begin{enumerate}
\item[(a)]
Because $V_D$ is a 3-disk and $(V_D)_{\gamma}$, obtained from $V_D$
by adding a 2-handle along $\gamma$, is a 3-disk with a hole, 
therefore adding this 2-handle does not change the KBSM (Lemma2.1(2)(ii)).
Therefore in ${\cal S}(V_D)$ and ${\cal S}(M - \{-1,1\} \times D)$
any relation of the form $L= sl_{\gamma}(L)$ holds. The embedding
$i: (M - \{-1,1\} \times D) \to M$ induces the homomorphism
of the KBSM\ $i_*: {\cal S}(M - \{-1,1\} \times D) \to {\cal S}(M)$,
thus $i_*(L) = i_*(sl_{\gamma}(L))$. Lemma 6.1 follows, as assumptions
of the lemma are chosen in such a way that any allowed sliding in $M$
is also a sliding in $M - (\{-1,1\} \times D)$.
\item[(b)]
$L \cap M_0$, being a  trivial knot, can be isotoped into $V_D$ without
changing an ambient isotopy type of the result of the sliding, which will
be also a trivial knot by part (a) of the lemma.
\item[(c)]
Any link $L$ in $M$ can be written in ${\cal S}(M)$
 as a linear combination of elements of $L^{gen}_{fr}$. A sliding 
described in Fig. 5.1 does not depend on the presentation of $L$
(or elements of $L^{gen}_{fr}$) so the lemma folllows. Notice
that the sliding relation of Fig. 5.1 performed on the link $D_L$
disjoint from $D$ holds already in ${\cal S}(M)$. 
\item[(d)] 
Let $D_L$ be a realization of a link $L$ in $M$.
As $D_L$ is arbitrary, one can assume that sliding has support
in $V_D$. Using relations from $I_0$ and  the conclusion of the 
part (c) of the lemma
together with Theorem 5.1 we can see that our slidings 
are performed on links in $M- \{-1,1\}\times D$ and have support in $V_D$,
so by (a) of  the lemma they do not introduce any new relation.\\
To visualize the assertion that modulo $I_0$ we need to slide only
links $D_L$ in $M- \{-1,1\}\times D$ with sliding support in $V_D$,
consider disks $D_1=\{-1\}\times D$ and $D_2=\{1\}\times D$ with 
$\gamma_i =\partial D_i$. Let $D_L$ be an arbitrary link in $M$ and
$sl_{\gamma}$ a sliding with support in $int (V_D)$. Slidings along
${\gamma}_1$ and ${\gamma}_2$ of the type described in Fig. 5.1, 
yield relations in ${\cal S}(M)$ satisfied in ${\cal S}(M)/I_0$
(as $\gamma_i$ is parallel to $\gamma$). These slidings allow us to
reduce $D_L$ to a linear combination of links (curve systems)
in $M-D_1 - D_2$. Furhermore the support of sliding 
$sl_{\gamma}$ is (unchanged) in $int (V_D)$.

\end{enumerate}
\end{proof}
As a corollary we get the main theorem for handlebodies.

\begin{corollary}[Main theorem for handlebodies]\label{6.2}\ \\
Let $D$ be a meridian disk of a handlebody $H_n$, ${\gamma}= \partial D$.
If $R$ is a ring with $(A^k-1)$ invertible for any $k$ then the embedding 
$j: (H_{n} - D) \to (H_n)_{\gamma}$ induces an isomorphism
$$j_*: {\cal S}(H_{n} - D) \to {\cal S}((H_n)_{\gamma}).$$
\end{corollary}

\begin{proof}
By the handle sliding lemma (Lemma 4.1) one has
${\cal S}((H_n)_{\gamma}) = 
{\cal S}(H_{n})/J  $ where $J$ is the submodule  of ${\cal S}(H_{n})$ 
generated by slidings $L-sl_{\gamma}(L)$.
By Lemma 6.1 $J=J_0$ for any generating set 
$L^{gen}_{fr}$ of ${\cal S}(H_{n})$. 
We can assume that $H_n = P_n \times [0,1]$ and $D = C \times [0,1]$ for 
an arc $C$ properly embedded in the disk with $n$ holes, $P_n$. 
Let $B(P_n)$ be the basis of 
${\cal S}(H_n)$ as described in Proposition 2.1(5). 
Let $B_i(P_n)$ be a subset of $B(P_n)$ composed of links with geometric
intersection number with $C$ equal to $i$. By Lemma 6.1, $J_0$ is
generated by sliding relations of Fig. 5.1, one relation for each
element of $B_i(P_n)$ for $i>0$. Because $B(P_n)$ is a basis of 
${\cal S}(H_n)$, therefore $B_0(P_n)$ is a basis of ${\cal S}(H_{n})/J_0$. 
On the other hand $B_0(P_n)= B(P_n-C)$ is also a basis of 
${\cal S}(H_{n} - D)$, thus $j_*$ is an isomorphism.
\end{proof}

We are ready now to prove the main theorem.\\
\begin{theorem}
Let $D$ be a meridian disk of $M$ and ${\gamma}= \partial D$.
If $R$ is a ring with $(A^k-1)$ invertible for any $k$ then the embedding 
$j: (M-D) \to M_{\gamma}$ induces an isomorphism 
$$j_*: {\cal S}(M- D) \to 
{\cal S}(M_{\gamma})$$ 
\end{theorem} 

\begin{proof}
$M_{\gamma}$ can be obtained from $(H_n)_{\gamma}$ by adding to
$M_{\gamma}$ some 2-handles (disjoint from the 2-handle added along
$\gamma$) and some 3-handles. By Lemma 4.1 and Proposition 2.1(2)
${\cal S}(M_{\gamma})$ is obtained from 
${\cal S}((H_n)_{\gamma})$ by sliding links generating 
${\cal S}((H_n)_{\gamma})$ along these 2-handles. 
Denote these slidings by $sl_h$. Consider any link $L$ in $(H_n)_{\gamma}$
and any sliding $sl_h$. We can choose a representative $D_L$ of $L$ so
that $D_L$ and $sl_h(D_L)$ are identical in the neighborhood of 
$S^2=D_{\gamma}$. By Lemma 6.1 we can present $D_L$ in 
${\cal S}((H_n)_{\gamma})$ as a linear combination of links which are 
disjoint from $S^2$ and differ from $D_L$ only in a small neighborhood 
of $S^2$. Thus the sliding relation $D_L - sl_h(D_L)$ is a linear combination
of sliding relations in $(H_n)_{\gamma} - S^2$. Therefore
$j_*: {\cal S}(M- D) \to {\cal S}(M_{\gamma})$ is an isomorphism.
The proof of Theorem 6.3 is complete.
\end{proof}
\begin{corollary}\label{6.4}
Let $R$ be a ring with $(A^k-1)$ invertible for any $k$ then
\begin{enumerate}
\item[(i)] If $S^2$ is a 2-sphere embedded in $M$ then the embedding
$j: M-S^2 \to M$ yields an isomorphism of the KBSM 
$j_*: {\cal S}(M- S^2) \to {\cal S}(M)$.
\item[(ii)] The embedding $M_1\# D^3 \sqcup M_2\# D^3 $ to $M_1\#M_2$
yields an isomorphism of the KBSM.
\item[(iii)] ${\cal S}(M_1\#M_2)$ is isomorphic to 
${\cal S}(M_1)\otimes {\cal S}(M_2)$.
\item[(iv)] ${\cal S}(S^1\times S^2) = R$
\end{enumerate}
\end{corollary}
\begin{proof}
The case (i) follows immediately from Theorem 6.3. The case (ii) 
corresponds to the case
of (i) when $S^2$ separates $M$. $M_1\#M_2 - S^2$ and $M_1\# D^3 \sqcup
M_2\# D^3 $ differ only by parts of their boundaries so 
their KBSM are the same. The case (iii) follows from (ii) by 
Proposition 2.1(3). The case (iv) follows from (i) as 
$S^1\times S^2 - S^2$ is a 3-disk with two holes. 
The case (iv) is also a special case of a general theorem in 
\cite{H-P-2}.
\end{proof}
\begin{corollary}\label{6.5}
If $S^2$ is a 2-sphere embedded in $M$ and the KBSM 
${\cal S}(M-S^2;Z[A^{\pm 1}],A)$ is free then the embedding
$j: M-S^2 \to M$ yields a monomorphism of the KBSM
$j_*: {\cal S}(M- S^2;Z[A^{\pm 1}],A) \to {\cal S}(M;Z[A^{\pm 1}],A)$.
\end{corollary} 
\begin{proof}
Let the set $\{x_{\alpha}\}$ be a  basis of ${\cal S}(M-S^2;Z[A^{\pm 1}],A)$. 
By the Universal Coefficient Property it is also a basis
of ${\cal S}(M-S^2;{\cal F}(A),A)$. By Theorem 6.3
the set $\{j_*(x_{\alpha})\}$ is a basis of ${\cal S}(M;{\cal F}(A),A)$. 
Therefore it is a $Z[A^{\pm 1}]$  linearly independent 
set in ${\cal S}(M;Z[A^{\pm 1}],A)$, and therefore $j_*$ is
a monomorphism.  
\end{proof}

\section{Generalizations and Speculations.}\label{7}
Theorem 1.2 does not hold for the ring $R=Z[A^{\pm 1}]$. As observed in
\cite{P-2}, Theorem 4.4, ${\cal S}(M_1\#M_2;Z[A^{\pm 1}],A)$
often contains a torsion part. The general description (generators and 
relators) of the KBSM is possible by Theorem 4.1 but to have a more 
meaningful description one should first analyze relative KBSM in
manifolds $M_1$ and $M_2$ and the KBSM of the connected sum would be
the sum of tensor products of ``reduced" relative skein modules. 
Finally one should be able to obtain for KBSM a Van Kampen-Seifert type 
theorem for 3-manifolds (glued along surfaces)
\footnote{The recent paper by W.Lofaro is a step in this direction
\cite{Lof}.}. The theorem 
could be reminiscent of Topological Quantum Field Theory 
formalism \cite{At}. We plan to give
the detailed description of the KBSM of  connected and disc sums of
3-manifolds in \cite{P-4}. Here we quote one, relatively simple result
(where there is no need to invoke the notion of the relative KBSM).
\begin{theorem}
$${\cal S}(H_n\# H_m) = {\cal S}(H_{n+m})/\cal I$$
where $\cal I$ is the ideal generated by expressions $z_k-A^6u(z_k)$, 
for any even $k\geq 2$, and $z_k \in B_k(P_{n+m})$, where 
$B_k(P_{n+m})$ is a subset of a basis $B(P_{n+m})$ composed of links 
with geometric intersection number with a disk $D$ separating
$H_n$ and $H_m$ equal to $k$. $u(z_k)$ is a modification of $z_k$
in the neighborhood of $D$, as shown in Fig. 7.1.
Our relation $z_k = A^6u(z_k)$, is a result of the sliding relation
$z_k = sl_{\partial D}(z_k)$ as illustrated  in Fig. 7.2.
\end{theorem}
   
\ \\ 
\centerline{\psfig{figure=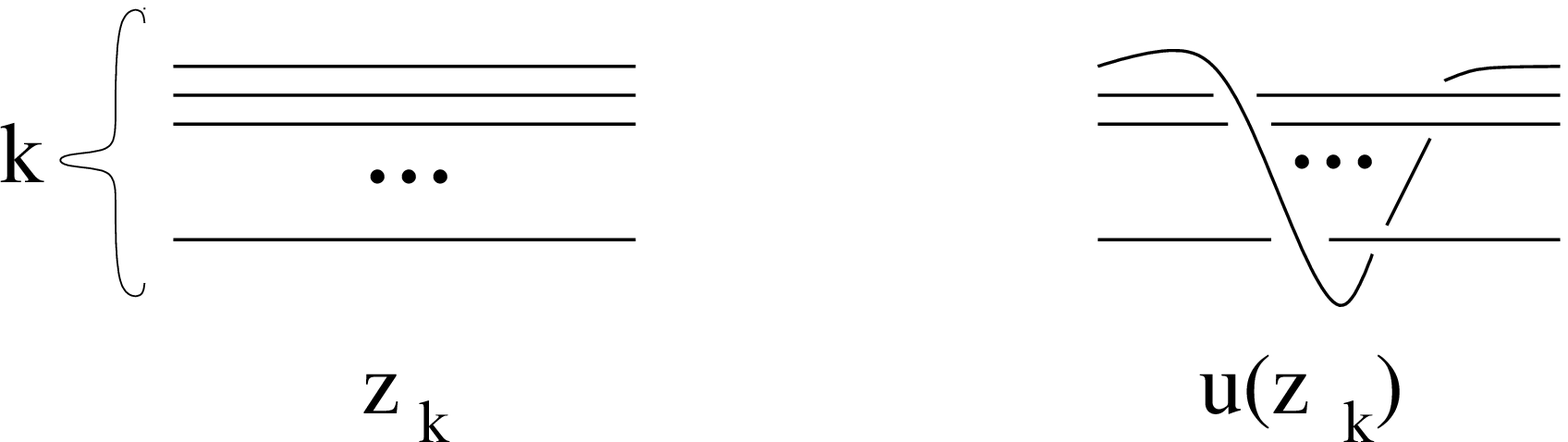,height=2.8cm}}
\begin{center} 
Fig. 7.1. 
\end{center} 
\ \\
\centerline{\psfig{figure=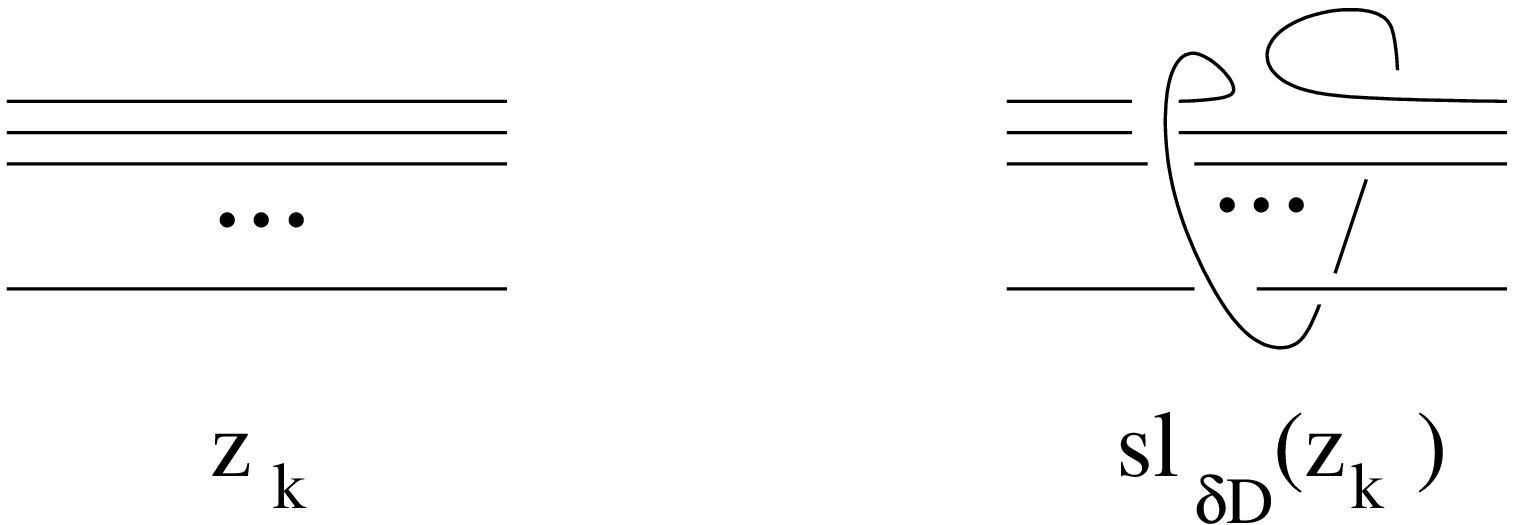,height=2.8cm}}
\begin{center} 
Fig. 7.2. 
\end{center}

\centerline{\it Department of Mathematics, University of Maryland}
\centerline{\it College Park, MD 20742}

The author is on leave from:\\
\centerline{\it Department of Mathematics, The George Washington University} 
\centerline{\it 2201 G Str. Funger Hall} 
\centerline{\it Washington, D.C. 20052} 
\centerline{\it przytyck@gwu.edu} 
\end{document}